\input amstex
\input amsppt.sty
\magnification=1200
\def\Aut{\operatorname{Aut}}

\def\phi{\varphi}
\def\alg{\operatorname{alg}}
\def\dist{\operatorname{dist}}
\def\Im{\operatorname{Im}}
\def\Re{\operatorname{Re}}
\define\pr{\operatorname{pr}}
\NoBlackBoxes

\topmatter
\title The Serre problem with Reinhardt fibers
\endtitle
\author Peter Pflug (Oldenburg) and W\l odzimierz Zwonek (Krak\'ow)
\endauthor
\abstract The Serre problem for a class of hyperbolic pseudoconvex
Reinhardt domains in $\Bbb C^2$ as fibers is solved.
\endabstract
\address
Carl von Ossietzky Universit\"at Oldenburg, Fachbereich
Mathematik, Postfach 2503, D-26111 Oldenburg, Germany
\endaddress
\address
Uniwersytet Jagiello\'nski, Instytut Matematyki, Reymonta 4,
30-059 Krak\'ow, Poland
\endaddress
\email pflug\@mathematik.uni-oldenburg, zwonek\@im.uj.edu.pl
\endemail
\rightheadtext{The Serre problem with Reinhardt fibers}
\thanks Research partially supported by the KBN grant No. 5 P03A 033 21
and\hfill \newline by the Nieders\"achsisches Ministerium f\"ur Wissenschaft
und Kultur, Az. 15.3 - 50 113(55) PL.
\endthanks

\endtopmatter
\document
Our aim is to discuss the Serre problem, i.e. the problem whether the
holomorphic fiber bundle $\pi:E\mapsto B$ with a Stein base $B$ and a Stein fiber
$F$ is Stein. For a comprehensive list of positive partial results to this problem see e.g. \cite{Siu}.

In our paper we consider this problem under the additional
assumption that the fiber $F$ is a pseudoconvex hyperbolic
Reinhardt domain in $\Bbb C^2$. Note that the first examples
showing that the answer to the Serre problem is in general
negative were constructed for Reinhardt fibers (see \cite{Sko},
\cite{Dem}, and \cite{Loeb}). Also first counterexamples with
bounded domains as fibers were found in the class of pseudoconvex
Reinhardt domains (see \cite{Coe-Loeb}).

We are interested in the problem, which bounded pseudoconvex Reinhardt
domains as fibers guarantee that the holomorphic fiber bundle with the
Stein basis is Stein, in other words for which bounded pseudoconvex Reinhardt domains
the answer to the Serre problem is positive.

Since in the class of pseudoconvex Reinhardt domains hyperbolicity (in the
sense of Carath\'eodory, Kobayashi or Brody) is equivalent to the
boundedness of domains (see \cite{Zwo~1}), it is natural that instead of
bounded we study the class of hyperbolic pseudoconvex Reinhardt domains.

Let us denote the class of Stein domains $D$
for which the answer to the Serre problem (with the fiber equal to $D$)
is positive by $\frak S$.

Now we may formulate our main theorem, which gives the characterization of hyperbolic pseudoconvex
Reinhardt
domains in $\Bbb C^2$ belonging to $\frak S$.

\proclaim{Theorem 1} Let $D$ be a hyperbolic pseudoconvex Reinhardt domain
in $\Bbb C^2$. Then $D\not\in\frak S$ if and only if $D$ is algebraically
equivalent to a Reinhardt domain $\tilde D\subset\Bbb C^2_*$ for which there
is a matrix $A\in\Bbb Z^{2\times 2}$ with the eigenvalues $\lambda$ and
$\frac{1}{\lambda}$, where $\lambda>1$, such that $$ \log\tilde
D=\{tv+sw:s>\phi(t),\;t>0\} \text{ {\rm(}or $\log\tilde
D=\{tv+sw:s>\phi(t),\;t<0\}${\rm)}}, $$ where $v,w\in\Bbb R^2$ are
eigenvectors corresponding to the eigenvalues $\lambda$ and
$\frac{1}{\lambda}$ and $\phi:(0,\infty)\mapsto[0,\infty)$
{\rm(}respectively, $\phi:(-\infty,0)\mapsto[0,\infty)${\rm)} is a convex
function satisfying the equality $\phi(t\lambda)=\frac{1}{\lambda}\phi(t)$,
$t\in(0,\infty)$ {\rm(}respectively, $t\in(-\infty,0)${\rm)}.
\endproclaim
Recall that the first known example of a bounded domain not
belonging to $\frak S$ was a domain from the class considered in Theorem 1.
More precisely, it was a domain associated to $A=\bmatrix 2 & 1\\
                                                          1 & 1  \endbmatrix$
and $\phi\equiv 0$ (defined on $(0,\infty)$) -- see
\cite{Coe-Loeb}. Later, D. Zaffran in \cite{Zaf} delivered other
domains not from $\frak S$ of the same type. Namely, he considered
domains associated to so-called 'even
Dloussky matrices' {i.e. $A=\bmatrix 0 & 1\\
                                 k_1 & 1 \endbmatrix\cdots\bmatrix 0 & 1\\
                                                                   k_{2s} & 1\endbmatrix$
                                                                   , $s, k_j\in\Bbb
                                                                   N\setminus\{0\}$},
                                                                    and with $\phi\equiv 0$.

In our considerations the key role in the proofs of positive results
(i.e. the facts that domains are from the class $\frak S$) will
be played by the criterion of Stehl\'e, which we formulate in the form
that we shall use in our paper.

\proclaim{Theorem 2 {\rm (see \cite{Ste} and \cite{Mok})}} Let $D$
be a domain in $\Bbb C^n$. If there exists a real-valued
plurisubharmonic exhaustion function $u$ on $D$ such that for any
$F\in\Aut D$ the function $u\circ F-u$ is bounded from above on
$D$, then $D\in\frak S$.
\endproclaim
$\Aut D$ denotes the group of holomorphic automorphisms of $D$.

Let us make a general remark. Below in the proofs we shall be
interested only in the cases when the group $\Aut D$ is not
compact; if $\Aut D$ is compact then $D\in\frak S$, which follows
from a general result (see \cite{K\"on} and \cite{Sib}).

The proof of Theorem 1 will be divided
into three different cases, depending on the number of
axis of $\Bbb C^2$ which intersect the domain $D$.

Formally, for a pseudoconvex Reinhardt domain $D\subset\Bbb C^n$ we define
$$
t:=t(D):=\#\{j\in\{1,\ldots,n\}:D\cap V_j\neq\emptyset\},
$$
where $V_j:=\{(z_1,\ldots,z_n)\in\Bbb C^n:z_j=0\}$, $j=1,\ldots,n$.

The three different cases we shall deal with in the proof of Theorem 1
correspond to the three possible values of $t$ (recall that $n=2$):

If $t=2$ (equivalently, $0\in D$) the result will simply follow from the
well-known sufficiency results for a domain to belong to $\frak S$. In fact
this case has already been done in \cite{K\"on}.

In the cases $t=1$ (these domains will always belong to $\frak S$)
and $t=0$ we shall concentrate on the structure of $\Aut D$. In
the case when $t=1$ there will only be three classes of model
domains for which the group is not compact (the result will follow
from \cite{Shi}). Two of the classes will be relatively simple to
deal with and the third class will consist of one special domain
for which we shall use the Stehl\'e criterion together with
Theorem 6.

In the case $t=0$ we shall use the result of \cite{Shi} to see that $\Aut D$
coincides with the group of algebraic automorphisms, $\Aut_{\alg}D$.
Studying the geometric structure of $D$ we shall see that there are two
classes of domains admitting non-compact automorphism groups. Because of the
geometry of the logarithmic image we call these two classes '{\it
parabolic}' and '{\it hyperbolic}'. In the hyperbolic case, which will
deliver us a negative answer to the Serre problem, we shall construct a
counterexample proceeding as in \cite{Coe-Loeb}. On the other hand the
parabolic case will be done similarly as the special case of the domain in
the case $t=1$.

\demo{Proof of Theorem 1 in case $t=2$} In this case $0\in D$ and
$D$ is bounded, so $D$ is Carath\'eodory complete (see \cite{Pfl})
and, consequently, because of \cite{Hir} $D\in\frak S$. As
mentioned earlier this case has already been done in \cite{K\"on}.
\qed
\enddemo

Before we go on to the proof of the case $t=1$ let us recall the description
of hyperbolic pseudoconvex Reinhardt domains, some results and notions
related to this class of domains, and the structure of automorphism groups
of such domains.

Recall that for a pseudoconvex Reinhardt domain in $\Bbb C^n$ the
logarithmic image of $D$ $$ \log D:=\{(x_1,\ldots, x_n)\in\Bbb
R^n:(e^{x_1},\ldots,e^{x_n})\in D\} $$ is convex.

\proclaim{Theorem 3 {\rm (see \cite{Zwo~1})}} Let $D$ be a pseudoconvex Reinhardt domain in $\Bbb
C^n$. Then the following conditions are equivalent:

-- $D$ is (Kobayashi, Carath\'eodory or Brody) hyperbolic,

-- $D$ is algebraically
equivalent to a bounded domain,

-- $\log D$ contains no straight lines and $D\cap V_j$ is empty or
hyperbolic {\rm(}as a domain in $\Bbb C^{n-1}${\rm)}, $j=1,\ldots,n$.
\endproclaim

\remark{Remark} Observe that the number $t$ remains fixed under algebraic
biholomorphism.
\endremark
\vskip 0.5cm
 In view of Theorem 3 we see that for a pseudoconvex Reinhardt
domain $D\subset\Bbb C_*^n$ $$ \text{ $D$ is hyperbolic if and only if $\log
D$ contains no straight lines.}\tag{1} $$ If $D$ is a pseudoconvex Reinhardt
domain in $\Bbb C^n$, then any element $\Phi\in\Aut_{\alg}D$ must be of the
following form $$ \Phi(z)=(b_1z^{A^1},\ldots,b_nz^{A^n}),\;z\in D, $$ where
$A=\bmatrix A^1\\ \cdot\\ \cdot\\ \cdot\\A^n\endbmatrix\in\Bbb Z^{n\times
n}$, $|\det A|=1$ and $b_1,\ldots,b_n\in\Bbb C_*$ (see \cite{Zwo~1}) -- here
for $\alpha=(\alpha_1,\ldots,\alpha_n)\in\Bbb Z^n$ we denote
$z^{\alpha}:=z_1^{\alpha_1}\cdot\ldots\cdot z_n^{\alpha_n}$ for $z\in\Bbb
C^n$ such that if $\alpha_j<0$ then $z_j\neq 0$.

Consequently, the mapping $\tilde\Phi(x):=Ax+\tilde b$, $x\in\log D$,
$\tilde b:=(\log|b_1|,\ldots,\log|b_n|)$, is an affine isomorphism of $\log
D$.

We may easily verify (from Cartan Theorem and Theorem 3) that for a
hyperbolic pseudoconvex Reinhardt domain $D$ the group $\Aut_{\alg} D$ is
not compact if and only if \vskip5mm

\item{(2)} there is a sequence $(\tilde\Phi_n)$ ($\tilde\Phi_n$
corresponds to $\Phi_n\in\Aut_{\alg}D$) such that for some (equivalently,
any) $x\in\log D$ we have $||\tilde\Phi_n(x)||\to\infty$ as $n\to\infty$.

\vskip5mm

Following the notation in \cite{Zwo~2}, for a pseudoconvex Reinhardt domain
$D\subset\Bbb C^n$ and $a\in\log D$ (chosen arbitrarily) we denote $$ \frak
C(D):=\{v\in\Bbb R^n:a+\Bbb R_+v\subset\log D\}. $$ Recall that $\frak C(D)$
is a closed convex cone with the origin at $0$, independent of $a$. It is
easy to verify that for any $\Phi\in\Aut_{\alg}D$, where $D$ is a
pseudoconvex Reinhardt domain in $\Bbb C^n$, $$ A(\frak C(D))=\frak
C(D),\tag{3} $$ where $A$ denotes the matrix associated to $\Phi$.

\remark{Remark} Consider a hyperbolic pseudoconvex Reinhardt domain $D$ in
$\Bbb C^2$ with $t=1$. We claim that in this case $\Aut_{\alg}D$ is compact.
Actually, take any $\Phi\in\Aut_{\alg}D$. Without loss of generality $D\cap
V_1=\emptyset$, $D\cap V_2\neq\emptyset$. Then one may easily verify from
the description of $\Aut_{\alg}D$ that $\Phi(z)=(b_1z_1^{\pm
1},b_2z_1^{\alpha_1}z_2)$, $z\in D$, for some $\alpha_1\in\Bbb Z$ and
$\Phi(\cdot,0)$ is a biholomorphism of $D\cap V_2$ (as a subdomain of $\Bbb
C$), from which we easily conclude that $\Aut_{\alg}D$ is compact.\endremark

The problem of characterization of automorphism groups of Reinhardt domains
was studied in \cite{Shi}; for Reinhardt domains with smooth boundary see
also \cite{Isa-Kra} and papers quoted there. The results obtained there
together with the above remarks lead us to the following description of
hyperbolic pseudoconvex Reinhardt domains in $\Bbb C^2$ with $t=1$ and
non-compact automorphism groups.

\proclaim{Theorem 4 {(\rm see \cite{Shi}, Theorem 5)}} Let $D$ be a
hyperbolic pseudoconvex Reinhardt domain in $\Bbb C^2$ with $t=1$. Then
$\Aut D$ is not compact if and only if $D$ is algebraically equivalent to
one of the domains: $$ \gather \{z_1\in\Bbb C:|z_1|<1\}\times\{z_2\in\Bbb
C:r<|z_2|<1\}=:\triangle\times P(r,1), \text{ where $0\leq r<1$,}\tag{4}\\
\{(z_1,z_2)\in\Bbb C^2:|z_1|<1,\;0<|z_2|<(1-|z_1|^2)^{p/2}\}\,\text{ for
some $p>0$,}\tag{5}\\ \{(z_1,z_2)\in\Bbb
C^2:0<|z_2|<\exp(-|z_1|^2)\}.\tag{6}
\endgather
$$
Moreover, when $D$ is as in \thetag{4}, then the group of automorphisms consists of
the mappings of the form
$$
D\owns(z_1,z_2)\mapsto(a(z_1),b(z_2))\in D,
$$ where $a$ is an automorphism of $\triangle$ and $b$ is the automorphism
of $P(r,1)$.

When $D$ is as in \thetag{5}, then the automorphism group consists of
the mappings of the form
$$
D\owns(z_1,z_2)\mapsto\left(\alpha\frac{z_1-\beta}{1-\bar\beta
z_1},\gamma\frac{(1-|\beta|^2)^{\frac{p}{2}}}{(1-\bar\beta z_1)^p}z_2\right)\in
D,
$$
where $|\alpha|=|\gamma|=1$, $|\beta|<1$.

When $D$ is as in \thetag{6}, then the automorphism group consists
of the mappings of the form $$ D\owns(z_1,z_2)\mapsto(\alpha
z_1+\beta,\gamma\exp(-2\alpha\bar\beta z_1-|\beta|^2)z_2)\in D, $$
where $|\alpha|=|\gamma|=1, \beta\in\Bbb C$.
\endproclaim

Let us formulate one more auxiliary result. In the proof of Lemma
6 (and later in the proof of Theorem 8) the important role will be
played once more by a result of S. Shimizu (which is combined
below with Theorem 3).

\proclaim{Theorem 5 {\rm (see \cite{Shi})}} Let $D$ be a
hyperbolic pseudoconvex Reinhardt domain in $\Bbb C_*^n$. Then
$\Aut D=\Aut_{\alg}D$.
\endproclaim

\proclaim{Lemma 6} Let $\tilde D$ be a pseudoconvex Reinhardt
domain in $\Bbb C^2$ such that
$$
\tilde D\subset\Bbb C\times R\cdot\triangle,\; \text{ for some }
R>0,\;\Bbb C_*\times\{0\}\subset\tilde D\text{ and }
(1,0)\not\in\frak C(\tilde D).\tag{7} $$ Then for any
$\Phi\in\Aut\tilde D$
$$
\Phi(\tilde D\cap(\Bbb C\times\{0\}))=\tilde D\cap(\Bbb
C\times\{0\}).
$$
Put $D:=\tilde D\setminus(\Bbb C\times\{0\})$. Assume additionally
that
$$
(0,0)\in\tilde D\text{ or }\frak C(D)=\Bbb R_+(0,-1).\tag{8}
$$
Then $\Aut D=\Aut \tilde D_{|D}$.
\endproclaim
\demo{Proof of Lemma 6} We prove the first part of the lemma. It
is sufficient to show the inclusion '$\subset$'. First we claim
that for any different points $(z_1,0),(\tilde z_1,0)\in\tilde
D\cap(\Bbb C\times\{0\})$ the equality
$\Phi_2(z_1,0)=\Phi_2(\tilde z_1,0)$ holds. Suppose the contrary.
Then it follows from \thetag{7} that for some two different points
$(z_1,0),(\tilde z_1,0)\in\tilde D\cap(\Bbb C\times\{0\})$ $$
0=k_{\Bbb C\setminus\{0\}}(z_1,\tilde z_1)\geq k_{\tilde
D}((z_1,0),(\tilde z_1,0))\geq
k_{R\cdot\triangle}(\Phi_2(z_1,0),\Phi_2(\tilde z_1,0))>0, $$
where $k_{\Omega}$ denotes the Kobayashi pseudodistance of
$\Omega$ -- contradiction.

Therefore, $\Phi_2(\tilde D\cap(\Bbb C\times\{0\}))=\{z_2^0\}$ for some
$|z_2^0|<R$. It is sufficient to show that $z_2^0=0$. Suppose the contrary.
Then the fact that $(1,0)\not\in\frak C(D)$ implies that the well-defined
holomorphic function $\Bbb C_*\owns z_1\mapsto\Phi_1(z_1,0)\in\Bbb C$ is
bounded, so constant. Therefore, $\Phi$ is constant on $\tilde D\cap(\Bbb
C\times \{0\})$ -- contradiction.

Assume now additionally \thetag{8}.  It follows from the first
part of the lemma that $(\Aut \tilde D)_{|D}\subset\Aut D$. Assume
for a while that each $\Phi\in\Aut D$ extends holomorphically onto
$\tilde D$. We shall prove that such an extension maps $\tilde D$
to $\tilde D$. Let $\Phi\in\Aut D$ and let $\tilde\Phi$ denote the
extension of $\Phi$ to $\tilde D$, $\tilde\Phi:\tilde
D\mapsto\overline{\tilde D}$.  Our aim is to show that
$\tilde\Phi(\tilde D)\subset\tilde D$. Note that in the case
$(0,0)\in\tilde D$ the existence of plurisubharmonic peak
functions for $\partial\tilde D$ together with the maximum
principle for subharmonic functions easily shows that
$\partial\tilde D\cap\tilde\Phi(\tilde D)=\emptyset$, which
finishes the proof in this case. So assume that $\frak C(D)=\Bbb
R_+(0,-1)$. Then $(0,0)\in\partial \tilde D$. Similarly as in the
previous case (use the plurisubharmonic peak functions and the
maximum principle for subharmonic functions) we see that $\tilde
\Phi(\tilde D)\cap(\partial\tilde D\setminus\{(0,0)\})=\emptyset$.
Suppose that $(0,0)\in\tilde\Phi(\tilde D)$. Then certainly
$(0,0)=\tilde\Phi(z_1^0,0)$ for some $z_1^0\in\Bbb C_*$. Let
$z_2^0\in\Bbb C_*$ be such that $(z_1^0,z_2^0)\in D$. Since
$\Phi\in\Aut D=\Aut_{alg}D$ (use Theorem 5) and because of the
equality $A(\frak C(D))=\frak C(D)$, where $A$ is the matrix
corresponding to $\Phi$, we get that $\tilde
\Phi(z_1^0,|z_2^0|\triangle)=(\tilde z_1^0,|\tilde
z_2^0|\triangle)$ for some $(\tilde z_1^0,\tilde z_2^0)\in D$,
which contradicts the continuity of $\tilde\Phi$ and the equality
$\tilde\Phi(z_1^0,0)=(0,0)$.

Therefore, to prove the other inclusion it suffices to show that
each $\Phi\in\Aut D$ extends holomorphically onto $\tilde D$. Let
$\Phi\in\Aut D$. Note that $\Phi_2$ is bounded, so it extends
holomorphically onto $\tilde D$. Therefore, we may expand $\Phi_2$
into the Hartogs-Taylor series in $\tilde D$: $$
\Phi_2(z_1,z_2)=z_2^{j_0}\sum_{j\geq j_0}c_j(z_1)z_2^{j-j_0}, $$
where $j_0\geq 0$ and $c_{j_0}\not\equiv 0$.

Write the Hartogs-Laurent expansion of $\Phi_1$ in $D$: $$
\Phi_1(z_1,z_2)=\sum_{j\in\Bbb Z}d_j(z_1)z_2^j. $$ Since $\Phi_1$ is not
constant, there is a $j\in\Bbb Z$ such that $d_j\not\equiv 0$. Note that
there is a $j\in\Bbb Z$ such that $d_k\equiv 0$ for any $k<j$. Actually,
otherwise the function $\Phi_1\cdot\Phi_2$ would be unbounded on $D$, which
would contradict \thetag{8}. Let $j_1$ denote the smallest $j$ satisfying
this property. To finish the proof it is sufficient to show that $j_1\geq
0$. Suppose the contrary. Then there is a $k\in\Bbb N$ such that
$kj_1+j_0<0$. But this implies that the function $\Phi_1^k\cdot\Phi_2$ is
unbounded on $D$, which contradicts \thetag{8}. \qed
\enddemo

Let us formulate a result we shall need in the proof of Theorem 1.

\proclaim{Theorem 7} Let $\tilde D$ be a pseudoconvex domain in
$\Bbb C^n$ and let $M$ be a pure one-codimensional analytic subset
of $\tilde D$. Put $D:=\tilde D\setminus M$. Assume, additionally,
that $\Aut \tilde D_{|D}=\Aut D$. Then the fact that $\tilde
D\in\frak S$ implies that $D\in\frak S$.
\endproclaim
\demo{Proof of Theorem 7} Assume that $D\notin\frak S$. So there
is a holomorphic fiber bundle $\pi\:E\mapsto B$ with a Stein basis
$B$ and the Stein fiber $D$ but $E$ is not Stein. We find an open
covering $(U_j)_{j\in J}$ of $B$ and a system of biholomorphic
mappings $(\Phi_j)_{j\in J}$, $\Phi_j\:\pi^{-1}(U_j)\mapsto
U_j\times D$ with $\pr _1\circ\Phi_j=\pi|_{\pi^{-1}(U_j)}$. Put
for $j$, $j'$, $j\neq j'$, $U_j\cap U_{j'}\neq \varnothing$,
$$g_{j,j'}=(g_{j,j';1},g_{j,j';2})\:(U_j\cap U_{j'})\times
D\mapsto (U_j\cap U_{j'})\times D , \quad (x,y)\mapsto
(\Phi_j\circ\Phi_{j'}^{-1})(x,y),$$ i.e. $g_{j,j';2}(x,\cdot
)\in\Aut (D)$ for all $x\in U_j\cap U_{j'}$. By assumption,
$g_{j,j';2}(x,\cdot)$ is the restriction of a holomorphic
automorphism $\tilde g_{j,j';2}(x,\cdot)\in\Aut (\tilde D)$.
Define now $$\tilde g_{j,j'}\:(U_j\cap U_{j'})\times\tilde
D\mapsto (U_j\cap U_{j'})\times\tilde D, \quad \tilde
g_{j,j'}(x,y):=(x,\tilde g_{j,j';2}(x,y)).$$ This mapping is
bijective and for fixed $x\in U_j\cap U_{j'}$ holomorphic in
$\tilde D$. Moreover, for $y\in D$ the mapping is holomorphic as a
function of $x$. Using Hartogs' theorem it follows that $\tilde
g_{j,j'}$ is biholomorphic. Obviously, the cocycle conditions
remain to be true. Therefore, we have obtained a new holomorphic
fiber bundle $\tilde\pi\:\tilde E\mapsto B$ over $B$ with fiber
$\tilde D$. Moreover, we may assume that $E\subset\tilde E$ and
$\tilde\pi|_E=\pi$. Then, in virtue of the assumption, $\tilde E$
is Stein. Using that $D=\tilde D\setminus M$ with a pure
one-codimensional analytic set $M$ it follows that $\tilde
E\setminus E$ is a pure one-codimensional analytic subset of
$\tilde E$ and, therefore (cf. \cite{Doc-Gra}, page 99, Satz 1),
Stein; contradiction. \qed
\enddemo

\demo{Proof of Theorem 1 for $t=1$} As earlier announced we consider only cases when $\Aut D$ is not
compact. When $D$ is as in \thetag{4} then one may
easily verify that the following exhausting function
$$
u(z):=\max\{-\log(1-|z_1|^2),-\log\dist(z_2,\Bbb C\setminus P(r,1))\}
$$
satisfies the assumptions of Stehl\'e's criterion.

When $D$ is as in \thetag{5}, then one may easily verify that the
following exhausting function
$$
\max\{\log|z_2|-\frac{p}{2}\log(1-|z_1|^2),-\log|z_2|\},\;z\in D
$$
satisfies the assumptions of the criterion of Stehl\'e.

Now assume that $D$ is of the form as in \thetag{6}.

Denote
$$
\tilde D=\{(z_1,z_2)\in\Bbb C^2:|z_2|<\exp(-|z_1|^2)\}.
$$
Therefore, we get from Lemma 6 (note that $\frak C(D)=\frak
C(\tilde D)=\Bbb R_+(-1,0)+\Bbb R_+(0,-1)$)
$$
\Aut \tilde D_{|D}=\Aut D.\tag{9}
$$

Now we prove that $\tilde D\in\frak S$.
Elementary calculations show that for any $\Phi\in\Aut \tilde D$
$$
\tilde u(\Phi(z))=\tilde u(z), \;z\in\tilde D,\tag{10}
$$ where $\tilde
u(z):=\log|z_2|+|z_1|^2$, $z\in\tilde D$.

Define ($\log^+|\lambda|:=\max\{\log|\lambda|,0\}$, $\lambda\in\Bbb C$) $$
u(z):=\max\{\rho(\tilde u(z)),\log^+|z_1|\}, z\in\tilde D, $$ where
$\rho:[-\infty,0)\mapsto[0,\infty)$ is a continuous, $C^2$-smooth on
$(-\infty,0)$, convex and increasing function such that $\lim_{t\to
0^-}\rho(t)=\infty$ (e.g. $\rho(t):=\frac{-1}{t}$, $t<0$). Then it is
trivial to see that $u$ is exhausting for $\tilde D$. Calculating the Levi
form of $\rho\circ\tilde u$ we see that $\rho\circ\tilde u$ is
plurisubharmonic on $D$ (and consequently, because of the Riemann extension
theorem, on $\tilde D$). Moreover, it follows from \thetag{10} and the form
of $\Aut \tilde D$ that for any $\Phi\in\Aut\tilde D$ the function
$u\circ\Phi-u$ is bounded from above on $\tilde D$, which implies, in view
of the criterion of Stehl\'e, that $\tilde D\in\frak S$.

Then, because of \thetag{9} we may make use of Theorem 7 to see
that $D\in\frak S$, too. \qed
\enddemo

At the moment we are left with the case $t=0$. Before we go on to
the proof of this case we show some auxiliary results. More
precisely, we characterize all hyperbolic pseudoconvex Reinhardt
domains $D$ in $\Bbb C_*^2$ with non-compact automorphism groups.

Additionally, for our future needs we give some necessary conditions on the
form of automorphisms in one of the cases.

\proclaim{Theorem 8} Let $D$ be a hyperbolic pseudoconvex
Reinhardt domain in $\Bbb C^2_*$ {\rm (}i.e. $t=0${\rm)}. Then
$\Aut(D)$ is not compact if and only if $D$ is algebraically
biholomorphic to a Reinhardt domain $\tilde D$ in $\Bbb C^2_*$ of
one of the following two types:
$$
\text{there are a matrix $A\in\Bbb Z^{2\times 2}$ and a number
$\beta_2\neq 0$ such that}\tag{11} $$
$$ \gather \log\tilde
D=\{tv+sw:t\in\Bbb R,s>\psi(t)\}\text{ {\rm(}if $\beta_2>0${\rm)}
or }\\ \log \tilde D=\{tv+sw:t\in\Bbb R,s<\psi(t)\} \text{
{\rm(}if $\beta_2<0${\rm)}}
\endgather
$$ and $1$ is the only eigenvalue of $A$ with the eigenvector $w$ {\rm(}so
$Aw=w${\rm)}, $Av=v+w$ for some $v\in\Bbb R^2$, and $\psi:\Bbb
R\mapsto\Bbb R$ is a convex {\rm(}or concave in the second
case{\rm)} function satisfying the property
$\psi(t+\beta_2)=t+\psi(t)$, $t\in\Bbb R$.

$$\text{there is a matrix $A\in\Bbb Z^{2\times 2}$ with the eigenvalues
$\lambda$ and $\frac{1}{\lambda}$, $\lambda>1$, such that}\tag{12}
$$
$$ \log\tilde D=\{tv+sw:s>\phi(t),\;t>0\} \text{ or } \log\tilde
D=\{tv+sw:s>\phi(t),\;t<0\}, $$ where $v,w\in\Bbb R^2$ are
eigenvectors corresponding to the eigenvalues $\lambda$ and
$\frac{1}{\lambda}$ and $\phi:(0,\infty)\mapsto[0,\infty)$
{\rm(}or $\phi:(-\infty,0)\mapsto[0,\infty)${\rm)} is a convex
function satisfying the equality
$\phi(t\lambda)=\frac{1}{\lambda}\phi(t)$, $t\in(0,\infty)$
{\rm(}or $t\in(-\infty,0)${\rm)}.

\vskip 5mm Moreover, in the case \thetag{11} each automorphism
$\Phi$ must be such that: $$ \tilde\Phi(x)=\tilde Ax+\tilde{\tilde
b},$$ where $\tilde A\in\Bbb Z^{2\times 2}$, $|\det\tilde A|=1$,
$\tilde Aw=w$ and one of three possibilites holds:

-- $\tilde A=\Bbb I_2$, $\tilde{\tilde b}=0$;

--  there is some $\tilde v\in\Bbb R^2$ such that $\tilde A \tilde
v=w+\tilde v$,

-- the number $-1$ is the second eigenvalue of $\tilde A$
with the corresponding eigenvector equal to $\tilde v$.

Additionally, in all cases, if we denote $x=tv+sw$ and $\tilde\Phi(x)=\tilde
tv+\tilde sw$ then $\tilde s-\psi(\tilde t)=s-\psi(t)$, $t,\tilde t,s,\tilde
s\in\Bbb R$.
\endproclaim

\subheading{Remark 9} In fact the domains representing two
different cases in \thetag{11} are actually algebraically
equivalent (use the biholomorphism $z^{-\Bbb I_2}$, where $\Bbb
I_2$ denotes the unit matrix).

The examples of functions $\phi$ from \thetag{12} are the functions defined
as follows $\phi(t):=\frac{a}{t}$, $t>0$ (or $t<0$), where $a$ is some
fixed number, $a\geq 0$ (or $a\leq 0$).

The examples of functions $\psi$ from \thetag{11} are the functions
defined as follows $\psi(t):=\frac{t(t-\beta_2)}{2\beta_2}$,
$t\in\Bbb R$.

One of the examples of matrices satisfying \thetag{12} has already been
given in the remarks after Theorem 1. More generally, the examples of
matrices satisfying \thetag{12} may be of the following form $A=\bmatrix k &
1\\
                    k-1 & 1
                    \endbmatrix$, $k\in \Bbb N$, $k\geq 2$
or even Dloussky matrices as defined in \cite{Zaf}.

Let us note that as an example of a matrix $A$
satisfying \thetag{11} we may take the matrix $\bmatrix 1 & 0 \\
                                                   k & 1
                                                   \endbmatrix$,
where $k\in\Bbb Z\setminus\{0\}$. Then $w=(0,1)$, $v=(\frac{1}{k},0)$.

Because of the geometry of $\log D$ let us call the domains satisfying
\thetag{11} of '{\it parabolic}' type and those satisfying \thetag{12} of
'{\it hyperbolic}' type.

\vskip5mm

\comment
In the proof of Theorem 7 the key role will be played
once more by a result of S. Shimizu.

\proclaim{Theorem 9 {\rm (see \cite{Shi})}} Let $D$ be a hyperbolic
pseudoconvex Reinhardt domain in $\Bbb C_*^n$. Then $\Aut D=\Aut_{\alg}D$.
\endproclaim
\endcomment

Now let us go on to the proof of Theorem 8.

\demo{Proof of Theorem 8} Assume that $\Aut D$ is not compact.
Then in view of \thetag{1}, \thetag{2} and \thetag{3} $\frak C(D)$
equals
$$
\text{$\Bbb R_+{v^{\prime}}$, $v^{\prime}\neq 0$ or $\Bbb
R_+v^{\prime}+\Bbb R_+w^{\prime}$, where $v^{\prime},w^{\prime}$
are linearly independent.}\tag{13}
$$

Let $\Phi\in\Aut D=\Aut_{\alg}D$. Denote the corresponding mapping
$\tilde \Phi(x)=Ax+\tilde b$, $x\in\Bbb R^n$.

First we claim that $A$ must have a positive eigenvalue. We
consider two possibilities as given in \thetag{13}. Consider the
first case $\frak C(D)=\Bbb R_+v^{\prime}$. Then from the
invariance $A(\frak C(D))=\frak C(D)$ we easily get that
$v^{\prime}$ is an eigenvector with the positive eigenvalue. So
assume that $\frak C(D)=\Bbb R_+v^{\prime}+\Bbb R_+w^{\prime}$. We
use once more the equality $A(\frak C(D))=\frak C(D)$ to see that
two cases have to be discussed, namely:
$$ A(\Bbb
R_{>0}v^{\prime})=\Bbb R_{>0}v^{\prime},\;A(\Bbb
R_{>0}w^{\prime})=\Bbb R_{>0}w^{\prime} \text{ or } A(\Bbb
R_{>0}v^{\prime})=\Bbb R_{>0}w^{\prime},\;A(\Bbb
R_{>0}w^{\prime})=\Bbb R_{>0}v^{\prime}.$$ Note that in the first
case we are done. In the second one using a continuity argument we
easily find the existence of a $u\in\Bbb R_{>0}v^{\prime}+\Bbb
R_{>0}w^{\prime}$ such that $A(\Bbb R_{>0}u)=\Bbb R_{>0}$, which
also finishes the proof of our claim.

Remark that if $A=\Bbb I_2$ then, because of \thetag{1}, $\tilde
b=0$.

Note also that if $A$ has a negative eigenvalue different from $-1$ then
taking instead of $\Phi$ the automorphism $\Phi^2$ we see that $\Aut D$ has
an element with the associated matrix $A$ having two positive eigenvalues,
both different from $1$. Therefore, we see that if $\Aut D$ consists of more
elements than those associated to $A=\Bbb I_2$ (and then automatically,
$\tilde b=0$) then $\Aut D$ must contain an element of one of the following
forms.

\comment
 Therefore, we see that except for the possibilty $A=\Bbb
I_2$ (and then $\tilde b=0$) the following three cases for the
form of $A$ are possible.
\endcomment

$$
\gather
A \text{ has two eigenvalues; one of them equals $1$ and the other
$-1$},\tag{14}\\
A \text{ has only one eigenvalue equal to $1$ and $A\neq\Bbb I_2$},\tag{15}\\
A \text{ has two positive eigenvalues $\lambda$ and
$\frac{1}{\lambda}$, $\lambda>1$.}\tag{16}
\endgather
$$

We claim the following:

If $\Aut D$ is not compact then $\Aut D$ contains an element such
that the corresponding matrix $A$ satisfies \thetag{15} or
\thetag{16}.

Actually, suppose the contrary. Then from our considerations any
automorphism of $D$ must be such that $A=\Bbb I_2$ or $A$ is of
the form \thetag{14}. Then any automorphism $\Phi\in\Aut D$ (with
$A\neq\Bbb I_2$) must be of the following form
$\tilde\Phi(x)=\alpha_1w-\alpha_2v+\beta_1w+\beta_2v$, where
$x=\alpha_1w+\alpha_2v$, $v,w\in\Bbb R^2$ are linearly
independent, where $w,v$ are eigenvectors of $A$ corresponding to
the eigenvalues $1,-1$ and $\beta_1,\beta_2\in\Bbb R$. It is easy
to verify that $\tilde\Phi^{(2)}(x)=(\alpha_1+2\beta_1)w+\alpha_2
v=x+2\beta_1v$. Therefore, we conclude that $\beta_1=0$. Moreover,
one may easily see that there is at most one automorphism with one
of the eigenvalue equal to $1$ (and the eigenvector equal to $w$)
and the other eigenvalue equal to $-1$ (and the eigenvector equal
to $v$). But the group consisting of only one element of this form
(and the identity) is certainly not compact. Therefore, there must
be some other automorphism with the same pair of eigenvalues (and
with other eigenvectors). Let us call the matrices corresponding
to these automorphisms by $A_1$ and $A_2$. Note that
$A_1A_2\neq\Bbb I_2$ (otherwise, $A_1=A_2$, which is impossible).
Since the matrix $A_1A_2$ corresponds to some automorphism of $D$
and the determinant of $A_1A_2$ equals $1$, we easily arrive at
the contradiction to our assumptions. Namely, it follows from our
considerations that either both eigenvalues are equal to $-1$,
which is impossible or they are both not real, which is
impossible, either, or they have two different real eigenvalues
with absolute values different from $1$, which is excluded in this
case, either.

Therefore, we proved our claim. \comment the non-compactness of
$\Aut D$ implies that there is some automorphism $\Phi$ such that
the corresponding matrix $A$ satisfies \thetag{15} or
\thetag{16}.\endcomment

Let us make one more remark. If we choose $A\in\Bbb Z^{2\times 2}$
with $|\det A|=1$ and with all eigenvalues different from $1$ then moving the domain, if
necessary, we may assume that $\tilde b=0$. In fact, since
$\det(A-I)\neq 0$ there is a vector $x_0\in\Bbb R^2$ such that
$Ax_0+\tilde b=x_0$. Consequently, for any $x\in\Bbb R^2$ the following
equalities hold
$$
\Phi(x)=Ax+\tilde b=A(x-x_0)+Ax_0+\tilde b=A(x-x_0)+x_0.
$$
Therefore, moving the coordinate system, if necessary, we may
assume that
$$
\tilde\Phi(x)=Ax,\tag{17}
$$
where $A$ is as above.

In addition to the previous remark note that, when the only eigenvalue is
$1$ and $A$ is not the identity, then some simplification of the form of
$\tilde \Phi$ is also possible. Namely, then there are linearly independent
vectors $v,w$ such that $Aw=w$ and $Av=v+w$. Then $(A-\Bbb I_2)(\Bbb
R^2)=\Bbb Rw$. Write any element $x\in\Bbb R^2$ in the form
$x=\alpha_1w+\alpha_2v$, $\tilde b=\beta_1w+\beta_2v$. Then there is an
$x_0\in\Bbb R^2$ such that $$ Ax_0+\tilde b=x_0+\beta_2v. $$ Consequently,
$$ \tilde\Phi(x)=Ax+\tilde b=A(x-x_0)+Ax_0+\tilde b=A(x-x_0)+x_0+\beta_2v,
$$ which implies that moving the coordinate system, if necessary, we may
assume that $$ \tilde\Phi(x)=Ax+\beta_2v.\tag{18} $$

Assume now that there is an automorphism of $D$ such that the
associated matrix $A$ satisfies \thetag{15}. Then we may also
assume that $\tilde\Phi$ satisfies \thetag{18}. There are linearly
independent vectors $v,w$ such that $Aw=w$ and $Av=w+v$. Note that
$$
\tilde\Phi^{(k)}(x)=x+k(\alpha_2w+\beta_2v)+\frac{k(k-1)}{2}\beta_2w,\;x=\alpha_1w+\alpha_2v,\;
k\in\Bbb Z.\tag{19} $$ Now \thetag{1}, \thetag{19}, and the
convexity of $\log D$ imply that $\beta_2\neq 0$ and for any
$t\in\Bbb R$ there is (exactly one) $s:=s(t)\in\Bbb R$ such that
$sw+tv\in\partial\log D$. Moreover, the convexity of $\log D$
together with \thetag{19} implies that if $\beta_2>0$ then
$sw+tv\in\log D$ for any $s>s(t)$ and if $\beta_2<0$ then
$sw+tv\in\log D$ for any $s<s(t)$. Denote $\psi(t):=s(t)$. Then
because of the equality $\tilde\Phi(\partial(\log
D))=\partial(\log D)$, the property \thetag{19} (applied for
$k=1$), and the convexity of $\log D$, we get the convexity (or
concavity) of the function $\psi$ and the property
$\psi(t+\beta_2)=\psi(t)+t$, $t\in\Bbb R$, which gives us the form
as in \thetag{11}. Note also that the domain as in \thetag{11} has
a non-compact automorphism group. In fact, note that
$\tilde\Phi(\log D)=\log D$ and $||\tilde\Phi^{(k)}(x)||\to\infty$
as $k\to\infty$, $x\in\log D$, where $\tilde\Phi(x)=Ax+\beta_2v$,
$x\in\log D$.

Consider now the case when there is an automorphism $\Phi$ of $D$ as in \thetag{16}.
Then because of \thetag{17} we may assume that $\tilde b=0$.

Consider any point $x=sw+tv\in\log D$, $t,s\in\Bbb R$, where
$v,w\in\Bbb R^2$ are eigenvectors corresponding to eigenvalues
$\lambda,\frac{1}{\lambda}$. Then $Ax=\frac{s}{\lambda}w+t\lambda
v$ and, consequently, for any $k\in\Bbb Z$ $$
A^k(x)=\frac{s}{\lambda^k}w+t\lambda^kv.\tag{20} $$ Taking $-\log
D$ instead of $\log D$, if necessary, (which corresponds to the
mapping $z^{-\Bbb I_2}$), we may assume that there is a vector
$x_0=s_0w+t_0v\in\log D$, where $s_0>0$, $t_0\neq 0$. Assume that
$t_0>0$ (the case $t_0<0$ goes along the same lines). Then it
easily follows from \thetag{20}, the convexity of $\log D$, and
\thetag{1} that $\log D\subset\{sw+tv:t,s>0\}$. Now one may easily
see from \thetag{20} and the convexity of $\log D$ that
$\{t>0:\text{ there is an $s>0$ such that $tv+sw\in\log D$}\}$ is
an open interval $(0,\infty)$. Moreover, for any $t>0$ there is
exactly one $s(t)\geq 0$ such that $sw+tv\in\log D$ for $s>s(t)$
and $sw+tv\not\in\log D$ for any $s<s(t)$. We define
$\phi(t):=s(t)$. The convexity of $\phi$ follows from the
convexity of $\log D$. The property
$\phi(t\lambda)=\frac{1}{\lambda}\phi(\lambda)$ easily follows
from the property \thetag{20}.

If we assume that $\log D$ is of the form as in \thetag{12} then it
follows from the properties of $A$ that $||A^k(x)||\to\infty$ as
$k\to \infty$, $x\in\log D$, which gives non-compactness of $\Aut
D$.

Now let us go to the study of the necessary form of the
automorphisms of $D$ in the case \thetag{11}.

Since the cones $\frak C(D)$ in both cases \thetag{11} and \thetag{12} are
not linearly isomorphic we easily conclude from the considerations that led
us to the construction of the domain as in \thetag{11} that each of the
automorphisms must be of one of the forms as in \thetag{14} or \thetag{15}
or its corresponding matrix must be the identity. Therefore, to finish the
proof it is sufficient to verify the invariance condition. Note that $\tilde
\Phi(\partial\log D)=\partial\log D$. Then elementary calculations show that
the invariance condition holds for all the possible automorphisms.
\qed
\enddemo

\demo{Proof of Theorem 1 for $t=0$} As noted earlier it suffices
to consider only the cases of non-compact $\Aut D$. As proven in
Theorem 8 there are two possibilities. We consider the first
(hyperbolic) one. We show that if $D$ is such that \thetag{12} is
satisfied then $D\not\in\frak S$.

We proceed as in \cite{Coe-Loeb}, we even follow the notation from
that paper. We define
$$
V:=\{\zeta_1 v+\zeta_2 w:\zeta_1,\zeta_2\in\Bbb C,\Im\zeta_1>0 \text{ (respectively, $\Im\zeta_1<0$), }
\Im\zeta_2>\phi(\Im\zeta_1)\}.
$$
It is obvious that $V/\Bbb Z^2$ is biholomorphic to $D$. We put
$\Omega:=\Bbb C\times V$.

We define the group of automorphisms $G_{\Bbb Z}$ induced by $\Bbb Z\times\Bbb Z^2$ on $\Omega$
as follows. Let $(\zeta_0,b_0)\in\Bbb Z\times\Bbb Z^2$. Then
$$
\Omega\owns(\zeta,b)\mapsto(\zeta+\zeta_0,A^{\zeta_0}b+b_0)=(\zeta+\zeta_0,
\zeta_1\lambda^{\zeta_0}v+\zeta_2\frac{1}{\lambda^{\zeta_0}}w+b_0)\in\Omega,
$$
where$(\zeta,b)=(\zeta,\zeta_1v+\zeta_2w)\in\Omega$.

Note that the fact that the functions defined above leave the set $\Omega$ invariant
follows from the properties of $\phi$. Namely,
$$
\frac{1}{\lambda^{\zeta_0}}\Im\zeta_2>\frac{1}{\lambda^{\zeta_0}}\phi(\Im\zeta_1)=
\phi(\Im\lambda^{\zeta_0})\text{ if and only if  }
\Im(\zeta_2)>\phi(\Im\zeta_1),
$$
for any $\zeta_0\in\Bbb Z$, $(\zeta_1,\zeta_2)\in\Bbb C^2$ with
$\Im\zeta_1>0$ (respectively, with $\Im\zeta_1<0$).

Now, we define the desired holomorphic fiber bundle $E:=\Omega/G_{\Bbb Z}$,
which has $V/\Bbb Z^2$ as the fiber and $\Bbb C/\Bbb Z$ as the basis.

Below we show that there is no plurisubharmonic exhaustion function on $E$.
Suppose the contrary. Let $u$ denote a plurisubharmonic exhaustion function
on $E$.

First recall that there is a family $(f_R)_{R>1}$ of holomorphic functions
$\bar\triangle\mapsto \Bbb C$ satisfying the following properties:
$$
\gather
0<\Im f_R<\pi\text{ on $\triangle$},\\
\Re f_R(0)=0 ,\\
\lim_{R\to 1^+}\int_{0}^{2\pi}e^{\pm Re f_R(e^{i\theta})}d\theta=\infty.
\endgather
$$
One may define $f_R(\zeta):=\log i\frac{R+\zeta}{R-\zeta}$ -- see \cite{Coe-Loeb}.

Note that $\phi(t\lambda^k)=\frac{1}{\lambda^k}\phi(t)$, $k\in\Bbb Z$,
$t\in[1,\lambda]$ (respectively, $t\in[-\lambda,-1]$) and $\phi$ is
continuous. Therefore, there is a constant $a>0$ (respectively, $a<0$) such
that $\frac{a}{t}>\phi(t)$, $t>0$ (respectively, $t<0$).

Now for any $R>1$ we find functions $g_R$ and $h_R$ holomorphic on
$\triangle$, continuous on $\bar\triangle$ and such that $$ \Im
g_R(\zeta)=e^{\Re f_R(\zeta)},\quad \Im h_R(\zeta)=ae^{-\Re
f_R(\zeta)},\;|\zeta|=1 $$ (respectively, $$ \Im g_R(\zeta)=-e^{\Re
f_R(\zeta)},\quad \Im h_R(\zeta)=-ae^{-\Re f_R(\zeta)},\;|\zeta|=1). $$
Since for any $|\zeta|=1$, the inequality $\phi(\Im g_R(\zeta))-\Im
h_R(\zeta)<0$ holds, the maximum principle for subharmonic functions (note
that $\phi\circ h$, where $h$ is harmonic on $\triangle$, is subharmonic on
$\triangle$) implies that for any $\zeta\in\bar\triangle$, the inequality
$\phi(\Im g_R(\zeta))-\Im h_R(\zeta)<0$ holds (or, equivalently,
$g_R(\zeta)v+h_R(\zeta)w\in V$, $\zeta\in\bar\triangle$).

Define
$$
\Psi_R(\zeta):=u([(\frac{f_R(\zeta)}{\log\lambda},g_R(\zeta)v+h_R(\zeta)w)]_{G_{\Bbb Z}}),\zeta\in
U_R\supset\bar\triangle.
$$
Certainly, $\Psi_R$ is subharmonic on some neighborhood of
$\bar\triangle$.

It follows from the definition of $E$ that $$
\Psi_R(\zeta)=u([(\frac{f_R(\zeta)}{\log\lambda}- \left[\frac{\Re
f_R(\zeta)}{\log\lambda}\right],\lambda^{-\left[\frac{\Re
f_R(\zeta)}{\log\lambda}\right]} g_R(\zeta)v+\lambda^{\left[\frac{\Re
f_R(\zeta)}{\log\lambda}\right]}h_R(\zeta)w)]_{G_{\Bbb Z}}), $$ for any
$\zeta\in\bar\triangle$ ($[x]$ in the inner brackets denotes the largest
integer not exceeding $x$).

Note that the real part of the first component in the formula above
is from the interval $[0,1)$ and its imaginary part is from the interval
$(0,\frac{\pi}{\log\lambda})$. Moreover, for $|\zeta|=1$
$$
\Im(\lambda^{-\left[\frac{\Re f_R(\zeta)}{\log\lambda}\right]}
g_R(\zeta))=e^{\Re
f_R(\zeta)-\log\lambda\left[\frac{\Re f_R(\zeta)}{\log\lambda}\right]}\in[1,\lambda)
$$
and, similarly,
$$
\Im(\lambda^{\left[\frac{\Re f_R(\zeta)}{\log\lambda}\right]}
h_R(\zeta))\in(\frac{a}{\lambda},a]
$$
(respectively,
$$
\gather
\Im(\lambda^{-\left[\frac{\Re f_R(\zeta)}{\log\lambda}\right]}
g_R(\zeta))\in(-\lambda,-1],\\
\Im(\lambda^{\left[\frac{\Re f_R(\zeta)}{\log\lambda}\right]}
h_R(\zeta))\in(\frac{-a}{\lambda},-a]).
\endgather
$$
Consequently, there is some constant $M\in\Bbb R$ such that for any $R>1$
$$
\Psi_R(\zeta)\leq M,\;|\zeta|=1.
$$

Then the maximum principle for subharmonic functions gives for any $R>1$
$$
\Psi_R(0)\leq M.
$$
Note that (the sign depends on one of two possible cases)
$$
\lim_{R\to 1^+}\Im g_R(0)=\pm\lim_{R\to
1^+}\frac{1}{2\pi}\int_{0}^{2\pi}e^{\Re
f_R(e^{i\theta})}d\theta=\pm\infty.
$$
Similarly,
$$
\lim_{R\to 1^+}\Im h_R(0)=\infty.
$$
Therefore, since $\Re f_R(0)=0$, $\Im f_R(0)\in(0,\pi)$,
and $u$ is exhaustive on $E$, we get
$\lim_{R\to 1^+}\Psi_R(0)=\infty$ -- contradiction. This finishes the
proof of the hyperbolic case.

\vskip5mm

Now we are left only with the parabolic case.

Applying the mapping $z^{-\Bbb I_2}$ we may reduce ourselves to the case
$\beta_2<0$ (and then $\psi$ is concave). Note that we may assume that $v$
and $w$ are from $\Bbb Q^2$ and the coordinates of $w$ are relatively prime.
Using some algebraic biholomorphism (mapping $w$ to $(-1,0)$ and being such
that the determinant of the corresponding matrix composing of integers is
one)  we may assume additionally that $w=(0,-1)$. Note that in this case
$\frak C(D)=\Bbb R_+(0,-1)$.

Consider now the domain $\tilde D:=\operatorname{int} \bar D$.
Note that $\tilde D\cap(\Bbb C\times\{0\})=\Bbb C_*\times\{0\}$,
$\tilde D\subset\Bbb C\times R\cdot\triangle$ for some $R>0$ and
$\frak C(\tilde D)=\Bbb R_+(0,-1)$. Then Lemma 6 implies that
$\Aut D=\Aut\tilde D_{|D}$. Therefore, as earlier, because of
Theorem 7, it is sufficient to show that $\tilde D\in\frak S$.

Now we define $$ u(z):=\max\{\rho(\tilde u(z)),\log|z_1|,-\log|z_1|\},\;z\in
\tilde D, $$ where $(\log|z_1|,\log|z_2|)=t(z)v+s(z)w$, $\tilde
u(z)=s(z)-\psi(t(z))$, $\rho:[-\infty,0)\mapsto\Bbb R$ is a continuous,
increasing and a convex function, $C^2$-smooth on $(-\infty,0)$ and
$\lim_{t\to 0^-}\rho(t)=\infty$ (e.g. $\rho(t)=\frac{-1}{t}$, $t<0$). Note
that assuming that $\psi$ is additionally $C^2$-smooth we may verify,
calculating the Levi form of $\rho\circ\tilde u$, that $\rho\circ\tilde u$
is plurisubharmonic on $D$. Then applying the standard approximation of a
concave function with the help of the increasing sequence of $C^2$-smooth
concave functions we get that $\rho\circ\tilde u$ is plurisubharmonic on $D$
without the additional assumption on its smoothness, too. Consequently, $u$
is plurisubharmonic on $D$, and then also on $\tilde D$.

It is clear that $u$ is an exhausting function for $\tilde D$. We claim that
for any $\Phi\in\Aut \tilde D$ $u\circ\Phi-u$ is bounded from above on
$\tilde D$.

Actually, take $\Phi\in\Aut\tilde D$. It follows from the
description of $\Aut D$ in Theorem~8 that $\rho(\tilde
u(\Phi(z)))=\rho(\tilde u(z))$, $z\in\tilde D$. One may also
verify that for any $\Phi\in\Aut \tilde D$
$\max\{\log|\Phi_1(z)|,-\log|\Phi_1(z)|\}-\max\{\log|z_1|,-\log|z_1|\}$
is bounded from above on $\tilde D$. Then the Stehl\'e criterion
applies and $\tilde D\in\frak S$. \qed
\enddemo

\subheading{Remark 10} As we saw in the proof of Theorem~1 there
were three non-trivial cases. The domain $\{(z_1,z_2)\in\Bbb
C^2:0<|z_2|<\exp(-|z_1|^2)\}$ for $t=1$ and the domains of
parabolic type for $t=0$ are domains for which the automorphism
group is non-compact. The proof that they belong to class $\frak
S$ relies upon the proof of belonging to the class $\frak S$ of
some larger domain. On the other hand the domains of hyperbolic
type are always not from $\frak S$ and the proof is based upon the
construction of Coeur\'e and Loeb.

It is natural to ask the question what happens in higher dimension. Is
there a similar geometric-like description of the class of hyperbolic
pseudoconvex Reinhardt domains from $\frak S$?

\enddocument

Consequently,
In addition to the previous remark note that when the only
eigenvalue is $1$ and $A$ is not the identity then some
simplification of the form of $\tilde \Phi$ is also possible.
Namely, then there are linearly independent vectors $v,w$ such
that $Aw=w$ and $Av=v+w$. Then $(A-\Bbb I_2)(\Bbb R^2)=\Bbb Rw$.
Write any element $x\in\Bbb R^2$ in the form
$x=\alpha_1w+\alpha_2v$, $\tilde b=\beta_1w+\beta_2v$. Then there
is an $x_0\in\Bbb R^2$ such that
$$
Ax_0+\tilde b=x_0+\beta_2v. $$ Consequently,
$$
\tilde\Phi(x)=Ax+\tilde b=A(x-x_0)+Ax_0+\tilde
b=A(x-x_0)+x_0+\beta_2v,
$$
which implies that moving the coordinate system, if necessary, we
may assume that
$$
\tilde\Phi(x)=Ax+\beta_2v.\tag{4}
$$

First, we prove that $\Aut D$ cannot contain any automorphism such
that the only eigenvalue of $A$ is $-1$ and $A\neq-\Bbb I_2$.

In fact, suppose the opposite. Then we assume that \thetag{3}
holds. There are linearly independent vectors $v,w$ such that
$Aw=-w$ and $Av=w-v$. Denote any point $x=\alpha_1w+\alpha_2v$.
Then one may verify that
$A^k(x)=(-1)^{|k|}x+(-1)^{|k|+1}k\alpha_2w$, $k\in\Bbb Z$,  which,
because of the convexity of $\log D$ leads to the contradiction
with \thetag{1}.

Assume now that there is an automorphism such that the associated
matrix $A$ satisfies the inequality $A\neq\Bbb I_2$ and the only
eigenvalue of $A$ is $1$. Then we may assume that $\tilde\Phi$
satisfies \thetag{4}. There are linearly independent vectors $v,w$
such that $Aw=w$ and $Av=w+v$. Note that
$$
\tilde\Phi^{(k)}(x)=x+k(\alpha_2w+\beta_2v)+\frac{k(k-1)}{2}\beta_2w,\;x=\alpha_1w+\alpha_2v,\;
k\in\Bbb Z.\tag{5}
$$
Now \thetag{1}, \thetag{5} and the convexity of $\log D$ imply
that $\beta_2\neq 0$ and for any $t\in\Bbb R$ there is (exactly
one) $s:=s(t)\in\Bbb R$ such that $sw+tv\in\partial\log D$.
Moreover, the convexity of $\log D$ together with \thetag{5} imply
that if $\beta_2>0$ then $sw+tv\in\log D$ for any $s>s(t)$ and if
$\beta_2<0$ then $sw+tv\in\log D$ for any $s<s(t)$. Denote
$\psi(t):=s(t)$. Then because of the equality
$\tilde\Phi(\partial(\log D))=\partial(\log D)$, the property
\thetag{5} (applied for $k=1$) and the convexity of $\log D$, we
get the convexity (concavity) of the function $\psi$ and the
property $\psi(t+\beta_2)=\psi(t)+t$, $t\in\Bbb R$. Note also that
the domain as defined in the theorem at the point \thetag{a} has
really non-compact automorphism group. In fact, note that $\log D$
is convex, $\tilde\Phi(\log D)=\log D$ and
$||\tilde\Phi^{(k)}(x)||\to\infty$ as $k\to\infty$, $x\in\log D$.

From now on we assume that for any $\Phi\in\Aut D$ the associated
matrix satisfies the equality $A=\Bbb I_2$ or $A=-\Bbb I_2$ or $A$
has two different eigenvalues.

\comment
Now we claim that there is a mapping $\Phi\in\Aut D$ such
that the matrix $A\in\Bbb Z^{2\times 2}$ from the description of
$\Phi$ satisfies $\det A=1$ and its eigenvalues are $\lambda$ and
$\frac{1}{\lambda}$, where $\lambda>0$, $\lambda\neq 1$.

First we show that $\Aut D$ contains no automorphism such that the
matrix associated to the automorphism has only one eigenvalue of
the second order such that the linear subspace corresponding to
this eigenvalue is one-dimensional. The only such a possibility is
that the eigenvalue equals $1$ (or $-1$) and then it easily
follows from the properties of such subspaces that
$\tilde\Phi^{(k)}(x)=x+k\alpha_2w+k\beta$ (or
$\tilde\Phi^{2k)}(x)=x+2k\alpha_2w$), where $x=\alpha_1v+\alpha_2
w$, for some choice of linear vectors $v,w\in\Bbb R^2$,
$\beta\in\Bbb R^2$ is fixed. We easily get contradiction with lack
of straight lines in $\log D$. From now on we assume that either
the matrix considered has two different eigenvalues or it has one
eigenvalue but there are two linearly independent eigenvectors
corresponding to it.
\endcomment

Additionally, assume that $\Aut D$ is not compact.

First we prove that

\vskip5mm
\item{(6)} $\Aut D$ contains no
automorphism such that $A=-\Bbb I_2$ and $\Aut D$ contains no
automorphism such that one of the eigenvalues of $A$ (and,
consequently, both) is not real. \vskip5mm

In fact, suppose the opposite. Let us find all possibilities of
automorphisms with non-real eigenvalues of the matrix associated.
Note that calculating the characteristic polynomial of $A$ we get
that the only possibilities of pairs of non-real eigenvalues are:
$\pm(\frac{k+i\sqrt{4-k^2}}{2}),\frac{k-i\sqrt{4-k^2}}{2})$, where
$k=a_{11}+a_{22}\in\{0,-1,1\}$. In any of the above cases there is
a $j\in\{2,3\}$ such that $A^j=-\Bbb I_2$. Note that if $A$ is the
matrix corresponding to the automorphism $\Phi$ then $A^j$ is the
matrix corresponding to the automorphism $\Phi^{(j)}$. Therefore,
if \thetag{6} did not hold then there would be an automorphism
$\Phi\in\Aut D$ such that the corresponding matrix equals $-\Bbb
I_2$. Making use of \thetag{3} applied to this automorphism we may
assume (moving the coordinate system, if necessary) that $\tilde
b=0$. But in view of \thetag{2} there is a sequence of points
$(x_n)\subset\log D$ with $||x_n||\to\infty$. We get that
$-x_n\in\log D$. This leads easily, because of convexity, to the
contradiction with \thetag{1}. Therefore, \thetag{6} is proved.

\bigskip

It is also trivial to see that $\Aut D$ contains no automorphism
of the form $\tilde\Phi(x)=x+\tilde b$, where $\tilde b\neq 0$
(use once more \thetag{1}).

\bigskip

Now we claim that if $\Aut D$ contains only automorphisms with
matrices $A$ having the eigenvalues equal to $1$ and $-1$ (and the
identity) then it cannot be compact. Actually, suppose the
opposite. Then for any automorphism $\Phi$ it must be of the
following form
$\tilde\Phi(x)=\alpha_1w-\alpha_2v+\beta_1w+\beta_2v$, where
$x=\alpha_1w+\alpha_2v$, $v,w\in\Bbb R^2$ are linearly
independent. It is easy to verify that
$\tilde\Phi^{(2)}(x)=(\alpha_1+2\beta_1)w+\alpha_2 v=x+2\beta_1v$.
Therefore, we get that $\beta_1=0$. Moreover, we easily get from
\thetag{1} that there is at most one automorphism with one of the
eigenvalue equal to $1$ (and the eigenvector equal to $w$) and the
other eigenvalue equal to $-1$ (and the eigenvector equal to $v$).
But the group consisting of only one element of this form (and the
identity) is certainly not compact. Therefore, there must be some
other automorphism with the same pair of eigenvalues (and with
other eigenvectors). Let us call the matrices corresponding to
these automorphisms by $A_1$ and $A_2$. Note that $A_1A_2\neq\Bbb
I_2$ (otherwise, $A_1=A_2$, which is impossible). Since the matrix
$A_1A_2$ corresponds to some automorphism of $D$ and the
determinant of $A_1A_2$ equals $1$, it follows from our
considerations that either both eigenvalues are equal to $-1$,
which is impossible or they are both not real, which is
impossible, either, or they have two different real eigenvalues
with absolute values different from $1$, which is excluded in this
case.

Therefore, the non-compactness of $\Aut D$ implies that if $D$ is
not as in \thetag{a} then there is an automorphism such that
$A\in\Bbb Z^{2\times 2}$, $|\det A|=1$ and $A$ has two real and
different eigenvalues both with the absolute value different from
$1$. Then because of \thetag{2} we may assume that $\tilde b=0$.

Taking $A^2$ instead of $A$ we may assume that both eigenvalues
are additionally positive, so that they are as in point
\thetag{b}.

Consider any point $x=sw+tv\in\log D$, $t,s\in\Bbb R$. Then
$Ax=\frac{s}{\lambda}w+t\lambda v$ and, consequently, for any
$k\in\Bbb Z$
$$
A^k(x)=\frac{s}{\lambda^k}w+t\lambda^kv.\tag{7}
$$
Taking $-\log D$ instead of $\log D$, if necessary, (which
corresponds to the mapping $\Bbb C_*^n\owns
z\mapsto(\frac{1}{z_1},\frac{1}{z_2})\in\Bbb C_*^2$), we may
assume that there is a vector $x_0=s_0w+t_0v\in\log D$, where
$s_0>0$, $t_0\neq 0$. Assume that $t_0>0$ (the case $t_0<0$ goes
along the same lines). Then it easily follows from \thetag{7}, the
convexity of $\log D$ and \thetag{1} that $\log
D\subset\{sw+tv:t,s>0\}$. Now one may easily see from \thetag{6}
and the convexity of $\log D$ that $\{t>0:\text{ there is an $s>0$
such that $tv+sw\in\log D$}\}$ is an open interval $(0,\infty)$.
Moreover, for any $t>0$ there is exactly one $s(t)\geq 0$ such
that $sw+tv\in\log D$ for $s>s(t)$ and $sw+tv\not\in\log D$ for
any $s<s(t)$. We define $\phi(t):=s(t)$. The convexity of $\phi$
follows from the convexity of $\log D$. The property
$\phi(t\lambda)=\frac{1}{\lambda}\phi(\lambda)$ follows easily
from the property \thetag{7}.

If we assume that $\log D$ is of the form as in \thetag{b} then it
follows from the property of $A$ that $||A^k(x)||\to\infty$ as
$k\to \infty$, $x\in\log D$, which gives non-compactness of $\Aut
D$. This completes the proof of Theorem.

\enddocument